\documentclass{amsart}

\usepackage{amsmath,amssymb}

\newtheorem{lemma}{Lemma}
\newtheorem{prop}[lemma]{Proposition}
\newtheorem{cor}[lemma]{Corollary}
\newtheorem{thm}[lemma]{Theorem}

\newtheorem{thm?}[lemma]{Theorem?}

\newtheorem{conj}[lemma]{Conjecture}

\newcommand{\F}{\mathbb{F}}

\newcommand{\ra}{\ensuremath{\rightarrow}}

\newcommand{\R}{\mathbb{R}}
\newcommand{\Z}{\mathbb{Z}}
\newcommand{\PP}{\mathbb{P}}
\newcommand{\C}{\mathbb{C}}
\newcommand{\Q}{\mathbb{Q}}

\newcommand{\Aut}{\operatorname{Aut}}

\newcommand{\OO}{\mathcal{O}}
\newcommand{\Ql}{\mathbb{Q}_{\ell}}
\newcommand{\Zl}{\mathbb{Z}_{\ell}}
\begin{document}
\title{Galois Groups Via Atkin-Lehner Twists}
\author{Pete L. Clark}
\address{1126 Burnside Hall \\ Department of Mathematics and Statistics \\
McGill University \\ 805 Sherbrooke West \\ Montreal, QC, Canada
H3A 2K6} \email{clark@math.mcgill.ca} \maketitle
\newcommand{\CT}{(\mathcal{C},\mathcal{T})}
\newcommand{\ttop}{\operatorname{top}}
\newcommand{\lcm}{\operatorname{lcm}}
\newcommand{\Div}{\operatorname{Div}}
\newcommand{\Gal}{\operatorname{Gal}}
\newcommand{\rank}{\operatorname{rank}}

\begin{abstract}
Using Serre's proposed complement to Shih's Theorem, we obtain
$PSL_2(\F_p)$ as a Galois group over $\Q$ for at least $614$ new
primes $p$.  Assuming that rational elliptic curves with odd
analytic rank have positive rank, we obtain Galois realizations
for $\frac{3}{8}$ of the primes that were not covered by previous
results; it would also suffice to assume a certain (plausible, and
perhaps tractable) conjecture concerning class numbers of
quadratic fields.  The key issue is to understand rational points
on Atkin-Lehner twists of $X_0(N)$.  In an appendix, we explore
the existence of local points on these curves.
\end{abstract}
\section{Introduction}
\noindent The notorious Inverse Galois Problem is to decide for
which finite groups $G$ there exists a Galois extension $L/\Q$
with $\Gal(L/\Q) \cong G$ (for short, ``$G$ occurs over $\Q$'').
The optimistic guess that every finite group occurs over $\Q$ is
natural for at least two reasons.  First, by a famous theorem of
Hilbert, it is enough to realize $G$ as the Galois group of a
regular extension $L/\Q(t)$.  Now for any field $K$ one says that
a finite group $G$ \emph{occurs regularly} over $K$ if it is the
Galois group of a regular extension $L/K(t)$, and there are many
fields -- e.g. $\C$, $\R$, $\Q_p\ldots$ -- over which every finite
group is known to occur regularly. \\ \indent There is also the
lure of inductive reasoning: it is known that many finite groups
-- e.g. solvable, symmetric and alternating groups -- occur as
Galois groups over $\Q$ (and, of course, no finite groups have
been shown \emph{not} to occur!). Still, the progress towards
realizing all groups has been anything but steady: some of the
``simplest'' simple groups are still not known to occur. Consider,
for instance, the family of groups $PSL_2(\F_p)$ as $p$ ranges
over prime numbers.
\\ \\
Over thirty years ago, Shih showed that $PSL_2(\F_p)$ occurs
regularly over $\Q$ if for some $N \in \{2,3,7\}$, the Kronecker
symbol $(\frac{N}{p})$ is equal to $-1$ \cite{Shih}. Later, Malle
showed that $(\frac{5}{p}) = -1$ is also sufficient for
$PSL_2(\F_p)$ to occur regularly over $\Q$ \cite{Malle}. Note that
these two results leave a density $\frac{1}{16}$ set of primes
unaccounted for.  To the best of my knowledge, no further
occurrences of $PSL_2(\F_p)$ over $\Q$, regular or otherwise, have
been established.
\\ \\
In 1988 Serre proposed a method of extending Shih's theorem to
realize new groups $PSL_2(\F_p)$ over $\Q$. This method is
discussed in his book \emph{Topics in Galois Theory} \cite{TGT}.
Referring to a calculation of Elkies, Serre remarks that the
method works to realize $PSL_2(\F_{47})$ over $\Q$; notice that $p
= 47$ is covered by Malle's result but not by Shih's. Strangely,
no additional examples of the method are given or asked for.
\\ \\
In this paper we analyze Serre's method and show that it works to
give realizations of $PSL_2(\F_p)$ over $\Q$ for many primes $p$
not obtainable by any previous result.
\section{Analysis of Serre's Method}
\noindent Serre's approach \cite[$\S 5.4$]{TGT} begins with the
following formulation of Shih's Theorem:
\begin{thm}(\cite[Theorem 8]{Shih})
\label{Shih} Let $p^* = (-1)^{\frac{p-1}{2}}p$, and $N \in \Z^+$
be such that $(\frac{N}{p}) = -1$.  Let $C(N,p)$ be the curve
obtained by twisting the modular curve $X_0(N)$ using the
Atkin-Lehner involution $w_N \in \Aut(X_0(N))$ and the quadratic
extension $Q(\sqrt{p^*})/\Q$.  Then there is a regular Galois
covering $Y \ra C(N,p)$, defined over $\Q$, with Galois group
$PSL_2(\F_p)$.
\end{thm}
\noindent In the case that $C(N,p) \cong \PP^1$, this means
precisely that $PSL_2(\F_p)$ occurs regularly over $\Q$. This
occurs for all $p$ when $N \in \{2, \ 3, \ 7\}$, and we recover
our earlier statement of Shih's Theorem.  More generally, when
$C(N,p)(\Q) \neq \emptyset$, there remains the possibility of
finding an irreducible specialization. Indeed, we have the
following
\begin{thm}(Serre)
\label{Serre}
With hypotheses as above, if $C(N,p)(\Q)$ is infinite, then there are infinitely many linearly disjoint Galois
extensions $L/\Q$ with Galois group $G \cong PSL_2(\F_p)$.
\end{thm}
\noindent
In other words, whereas Shih's theorem says that it suffices for $C(N,p) \cong \PP^1$, Theorem \ref{Serre} says
that it also suffices for $C(N,p)$ to be an elliptic curve of positive rank.
\\ \\
$X_0(N)$ has genus one for $N \in \{11, \ 14, \ 15, \ 17, \ 19, \
20, \ 21, \ 24, \ 27, \ 32, \ 36,  \ 49\}$, so in order to realize
$PSL_2(\F_p)$ for primes $p$ not covered by the results of Shih
and Malle, the values of $N$ to consider are $11, \ 17$ and $19$.
\begin{thm}
\label{GoodNews}
Let $N = 11$ or $19$.  For all primes $p$,
$C(N,p)(\Q) \neq \emptyset$, so $C(N,p)$ can be given the
structure of a rational elliptic curve.  More precisely, $C(N,p)$
is the quadratic twist of $X_0(N)$ by $p^*$.  It
follows that: \\
$\bullet$ For primes $p \equiv 1 \pmod 4$, $(\frac{N}{p}) = -1$ $\iff$ $C(N,p)$ has odd analytic rank. \\
$\bullet$ For primes $p \equiv -1 \pmod 4$, $(\frac{N}{p}) = -1$ $\iff$ $C(N,p)$ has even analytic rank.
\end{thm}
\noindent Proof: Let $\sigma$ denote the nontrivial element of
$\Gal(\Q(\sqrt{p^*})/\Q)$.  Then $C(N,p)(\Q)$ consists of those
points $P$ of $C(N,p)(\Q(\sqrt{p^*}))$ for which $w_N(\sigma(P)) =
P$.  In particular, any $P \in X_0(N)(\Q)$ which is a fixed point
of $w_N$ remains $\Q$-rational on $C(N,p)$.  For any squarefree
$N$, let $S_1$ (resp. $S_2$) be the set of $\C$-isomorphism
classes of elliptic curves with CM by the maximal order of
$\Q(\sqrt{-N})$ (resp. by $\Z[\sqrt{-N}]$); these sets are
distinct if and only if $-N \equiv 1 \pmod 4$.  It is not hard to
see that $S_1$ and $S_2$ each gives rise to a set of $w_N$-fixed
points -- the $\Gamma_0(N)$ structure is given by $E[\sqrt{-N}]$
--  and that if $N
> 3$, $S_1 \cup S_2$ gives all the $w_N$-fixed points (e.g.
\cite[Prop. 3]{Ogg}).  The set $S_1$ (resp. $S_2$) forms a
complete $\Gal(\overline{\Q}/\Q)$ orbit.  So there are
$\Q$-rational $w_N$-fixed points exactly when $\Q(\sqrt{-N})$ has
class number $1$, which is the case for $N = 11$ and $N = 19$ (but
not for $N = 17$).  \\
\indent Thus, for $N = 11$ or $19$, let $O$ be the unique fixed
point of $w_N$ that parameterizes an elliptic curve with
$\Z[\frac{1+\sqrt{-N}}{2}]$-CM, so that $(X_0(N),O)$ and
$(C(N,p),O)$ are rational elliptic curves.  Since the
$j$-invariant of $X_0(N)$ is neither $0$ nor $1728$, the group of
automorphisms of $(X_0(N),O)$ is $\pm 1$, and it follows that
under our identifications we simply have $w_N = -1$. In other
words, the twist of $X_0(N)$ via $w_N$ and $\Q(\sqrt{p^*})/\Q$ is
just the quadratic twist by $p^*$ in the usual sense.  Finally,
the sign of the functional equation for both $X_0(11)$ and
$X_0(19)$ is $+1$, so the sign of the functional equation for
$C(N,p)$ is $\chi_{p^*}(-N) = (\frac{-N}{p})$.  The result
follows.
\\ \\
We get as an immediate consequence our main result:
\begin{cor}
Assume that rational elliptic curves with odd analytic rank have positive Mordell-Weil rank.  Then for every
$p$ with $p \equiv 1 \pmod 4$ and which is a quadratic nonresidue \emph{either} modulo $11$ or modulo $19$,
$PSL_2(\F_p)$ occurs infinitely often as a Galois group over $\Q$.
\end{cor}
\noindent However, when $N = 17$ Serre's method fails for local
reasons:
\begin{prop}
\label{17}
For every prime $p$ such that $\left(\frac{17}{p}\right)
= -1$, $C(N,p)(\Q_{17}) = \emptyset$.
\end{prop}
\noindent Remark: The proof that follows is purely computational.
More recently we have found other approaches: see Corollary
\ref{sharp17} and Theorem \ref{AtN} in the Appendix.
\\ \\
Proof: For any squarefree $N$, the involution $w_N$ has at least
one fixed point, so that the quotient map $X_0(N) \ra X^+_0(N) :=
X_0(N)/w_N$ is always ramified.  In particular, when $X_0(N)$ has
genus one, $X^+_0(N)$ has genus zero, and in fact --  since
$X_0(N)$ always has $\Q$-rational cusps -- $X_0^+(N) \cong \PP^1$.
In particular, $w_{17}$ gives a hyperelliptic involution on
$X_0(17)$, and $\Q(X_0(17)) = \Q(x,y)$, where $y^2 = P(x)$ for
some quartic polynomial $P \in \Q[x]$ and $w: (x,y) \mapsto
(x,-y)$. Explicit polynomials have been computed by Elkies
(unpublished) and Gonz\'alez; by \cite[p. 794]{Gonzalez}, we may
take $P(x) = x^4+2x^3-39x^2-176x-212$. Thus $C(17,p)$ is given by
the equation
\[p^* y^2 = x^4+2x^3-39x^2-176x-212.\] The condition
$(\frac{17}{p}) = -1$ ensures that any two of the numbers $p^*$
differ (multiplicatively) by an element of $\Q_p^{\times 2}$,
hence all such curves $C(17,p)$ are $\Q_{17}$-isomorphic to a
single curve, say $C(17,5)$, and we are reduced to deciding
whether this particular hyperelliptic quartic curve has
$\Q_{17}$-rational points.  \\ \indent As a general principle, one
knows that all Diophantine problems over local fields are
decidable; in this case an analysis involving Hensel's Lemma shows
that (after ruling out rational points at infinity and
projectivizing) it suffices to study a corresponding congruence
modulo $17^5$. Thus in practice we will want computer assistance.
Luckily for us, the issue of local points on hyperelliptic
quartics arises in $2$-descent, so that a relatively efficient
algorithm for this -- first described in the early work of Birch
and Swinnerton-Dyer -- has been implemented in several elliptic
curve software packages. For instance, a query to John Cremona's
program {\tt ratpoints} results in the instantaneous response that
$C(17,5)$ fails to have points rational over $\Q_{17}$ (and indeed
also over $\Q_5$). This completes the proof.
\section{Examples}
\noindent Fix $N \in \{11, 19\}$.  If $p \equiv 1 \pmod 4$ is a
prime satisfying $(\frac{N}{p}) = -1$, then we have seen that the
analytic rank of $C(N,p)$ is odd, so that it is, to say the least,
widely believed that $C(N,p)$ has infinitely many $\Q$-rational
points.  In any given case one can, in principle, verify this just
by finding a rational point of infinite order.  Better yet:
whenever the analytic rank is equal to $1$, we know -- thanks to
the work of Gross-Zagier and Kolyvagin -- that the rank is equal
to one.  And it is easy to check that the analytic rank is $1$:
one needs only to check that the special value $L'(C(N,p),1)$ is
nonzero, which is amenable to approximate computation.  Moreover,
when the sign in the functional equation is $-1$, the prevailing
wisdom is that the analytic rank should be equal to $1$ ``most of
the time.''  \\ \indent Indeed, among all primes $p \leq 3 \times
10^5$ for which
\[p \equiv 1 \pmod 4, \ \left(\frac{2}{p}\right) = \left(\frac{3}{p} \right) =
\left(\frac{5}{p}\right) = \left(\frac{7}{p}\right) = 1,\] it is
\emph{never} the case that the analytic rank of either $C(11,p)$
or $C(19,p)$ is odd and greater than one.  This yields a list of
$612$ prime numbers -- the smallest is $p = 1009$ -- such that
$PSL_2(\F_p)$ newly occurs as a Galois group over $\Q$.
\\ \\
On the other hand, when $p \equiv -1 \pmod 4$, the analytic rank
of $C(N,p)$ is noticeably less averse to being \emph{even} and
greater than zero.\footnote{At least apparently: approximate
calculation can strongly suggest, but not prove, that an elliptic
curve has analytic rank $\geq 2$.} Elkies' computation of $\rank
C(11,47) = 2$ is an example of this.  Among $p < 1000$, $C(11,p)$
also has (apparent) analytic rank $2$ for $103, \ 599,$ and $683$.
The phenomenon is robust enough to persist upon enforcing the
congruence conditions $(\frac{2}{p}) = (\frac{3}{p}) =
(\frac{5}{p}) = (\frac{7}{p}) = 1$; for
\begin{equation}
p = 4079, \ 5591, \ 6719, \ 10391, \ 19319, \ 24359, \ 26759,
\end{equation}
either $C(11,p)$ or $C(19,p)$ has (apparent) analytic rank $2$.  \\ \\
Example 1:
\[[51362438166007626829703:\! -4948233782238353787199293:\! 5697234033382001683] \]
is a nontorsion point on the curve
\[C(11,4079): Y^2Z + YZ^2 = X^3 - X^2Z - 171928490XZ^2+1571689994520Z^3,\]
so $PSL_2(\F_{4079})$ occurs as a Galois group over $\Q$.
\\ \\
Example 2: $[-99184162:21162527913:10648] $
is a nontorsion point on the curve
\[C(19,5591): Y^2Z + YZ^2 = X^3 + X^2Z - 291753289XZ^2 + 2040511796399Z^3,\]
so $PSL_2(\F_{5591})$ occurs as a Galois group over $\Q$.
\\ \\
On the basis of these computations, it seems reasonable to
conjecture that $C(N,p)$ has positive rank for infinitely many
primes $p \equiv -1 \pmod 4$.
\section{Some Remarks}
\noindent There is a large literature on the variation of ranks of
elliptic curves in a family of quadratic twists, but comparatively
little has been said about the case of restricting to twists by
prime numbers.  Worth mentioning in this regard is the work of Ono
and Skinner \cite{OS}, which shows that for certain rational
elliptic curves a positive proportion of prime twists have rank
zero. There is as yet no analogous result for positive rank.  Note
however that Vatsal has shown that a positive proportion of all
quadratic twists of $X_0(19)$ have rank $1$ \cite{Vatsal}.  His
argument could be readily adapted to the case of twists by primes
in a given congruence class -- and hence to give, unconditionally,
a new positive density set of primes $p$ for which $PSL_2(\F_p)$
occurs over $\Q$ -- provided we knew the following:
\begin{conj}
Fix coprime positive integers $m$ and $M$, and let $F(X)$ be the number of primes satisfying $p \leq X$,
$p \equiv m \pmod M$, and such that the class number of $\Q(\sqrt{-3p})$ is indivisible by $3$.  Then $F(X) \gg X/\log X$.
\end{conj}
\noindent This conjecture may be within reach; by taking $M = 4$
and replacing $\Q(\sqrt{-3p})$ with $\Q(\sqrt{p})$, we get a
theorem of Belabas and Fouvry \cite{BF}.
\\ \\
\noindent Professor Shih has made me aware of the relevance
 of his later paper \cite{Shih2} to the
 present work.  Especially, this note can be viewed as responding
 to \cite[Remark 4]{Shih2}.  What does not figure in \cite{Shih2}
 is the dichotomy between $N = 11, 19$ and $N = 17$ arising from
 the fact that the curves $C(N,p)$ have $\Q$-rational points in
 the former case but not in the latter case.  Deciding
 which Atkin-Lehner twists of $X_0(N)$ have points rational over
 $\Q$ (or even over all of its completions) is an interesting Diophantine
 problem which we address in the
 Appendix of this paper.
\\ \\
The calculations of Section 3 have been extended, thanks to the
help of Nick Rogers: in particular, for all the primes $p$ in (1),
$\rank C(11,p) + \rank C(19,p) > 0$. Conversely, after a wider
search, we have still not been able to find any primes for which
either $C(11,p)$ or $C(19,p)$ has odd analytic rank $\geq 3$.  It
has been suggested to me by several people that it is
``implausible'' for, e.g., the ranks to be bounded in a family of
quadratic twists by \emph{prime numbers}. In response to this, I
would like to say that I do not necessarily agree: although I am
too far from being sufficiently expert in this matter to have an
opinion of my own, I was not able to find in the literature any
conjecture, heuristic or model (let alone any theorem) which would
rule out, e.g., the statement that $\rank C(19,p) = 1$ for all
primes $p \equiv 1 \pmod 4$.  In any case, the phenomenon seems to
be worthy of further investigation.
\section*{Appendix: Local Points on $C(N,p)$}
\noindent We suppose that $N$ is squarefree and $p$ is a prime
number such that $(\frac{N}{p}) = -1$.
\\ \\
The cusps of $X_0(N)$ are $\Q$-rational but freely permuted by
$w_N$, so that if $\Q(\sqrt{-N})$ has class number greater than
one there are no ``obvious'' $\Q$-rational points on $C(N,p)$.
Indeed, Shih found that in certain cases $C(N,p)$ fails to have
points even over certain completions of $\Q$.
\\ \\
For any nonsquare integer $m$, let $C(N,m)$ denote the twist of
$X_0(N)$ by $w_N$ and $\Q(\sqrt{m})/\Q$.  As Ellenberg has pointed
out \cite[Problem A]{Ellenberg}, it is an interesting general
problem to classify the ``deficient places'' of $C(N,m)$, i.e.,
the primes $\ell \leq \infty$ such that $C(N,m)(\Q_l) =
\emptyset$.  We shall give some results in this direction.

\begin{prop}
\label{real} $C(N,p)(\R) \neq \emptyset$.
\end{prop}
\noindent Proof: Let $\OO$ be the ring of integers of
$\Q(\sqrt{-N})$.  Then $E := \C/\OO$ is an elliptic curve over the
complex numbers such that $P = (E,E[\sqrt{-N}])$ gives a
$w_N$-fixed point.  But complex conjugation on $\C$ induces an
antiholomorphic involution on $E$, and it follows easily that
$j(E) \in \R$ and that $P \in X_0(N)(\R)$.  Thus every
Atkin-Lehner twist of $X_0(N)$ has $\R$-rational points.
\\ \\
\noindent Recall that a (nonsingular, projective) curve $C_{/\Ql}$
is said to have \emph{good reduction} if there exists a smooth
arithmetic surface $\mathcal{C}_{/\Zl}$ with generic fiber
isomorphic to $C$.  If $C$ has positive genus, then it admits a
unique minimal model $\mathcal{C}_{/\Zl}$, whose smoothness is
equivalent to the good reduction of $C$.  On the other hand, there
are only two curves of genus zero over $\Ql$: $\PP^1$ (which has
good reduction), and the twisted form corresponding to the unique
division quaternion algebra over $\Q_{\ell}$ (which does not have
good reduction).
\begin{prop}
\label{smooth} Fix $\ell$ prime to $Np$.  Then $C(N,p)$ has good
reduction over $\Ql$.  In particular, if $X_0(N)$ has genus at
most one, $C(N,p)(\Ql) \neq \emptyset$.
\end{prop} \noindent Proof: When
$X_0(N)$ has genus zero -- i.e., when $N \in \{2, \ 3, \ 5, \ 6, \
7, \ 10, \ 13 \}$ -- work of Shih gives a more precise result.
Indeed, the genus zero curve $C(N,p)_{/\Q}$ is classified up to
isomorphism by a quaternion algebra.  By \cite[Prop. 10]{Shih},
this quaternion algebra is given by the Hilbert symbol $\langle
c_N, p^* \rangle$, where
\[c_2 = c_3 = 1, \ c_5 = 125, \ c_6 = 18, \ c_7 = 49, \ c_{10} = 5,
 \ c_{13} = 13. \]
The genus zero case follows (after a small calculation, when $\ell
= 2$).
\\ \indent
Let us now assume that $X_0(N)$ has positive genus.  By a well
known theorem of Igusa, $X_0(N)$ has good reduction over
$\Z_{\ell}$, so $C(N,p)_{/\Q(\sqrt{p^*})}$ has good reduction at
the places over $\ell$.  Moreover $\ell$ is unramified in
$\Q(\sqrt{p^*})$, so that $C(N,p)_{/\Q_{\ell}}$ has good reduction
after an unramified base change.  Since formation of the minimal
model commutes with unramified base change and smoothness can be
checked on geometric fibers, this means that the minimal model
$\mathcal{C}(N,p)_{/\Z_{\ell}}$ is itself smooth.  \\ \indent
Finally, if $X_0(N)$ has genus one, then, since every smooth curve
of genus at most one over the finite field $\F_{\ell}$ has an
$\ell$-rational point (e.g., by the Weil bounds), $C(N,p)(\Ql)
\neq \emptyset$ by Hensel's Lemma.
\\ \\
The next result is a generalization of \cite[Prop. 10]{Shih}.
\begin{thm}(Gonz\'alez)
\label{obstruction} $C(N,p)_{/\Q}$ admits a finite morphism to the
genus zero curve with corresponding quaternion algebra $\langle
c_N , \ p^* \rangle$, where $c_N =  N^{\frac{12}{\gcd(12,N-1)}}$.
\end{thm}
\noindent Proof: See \cite[Thm. 6.2]{Quer}. The argument is
analytic in nature: one constructs a $\Gamma_0(N)$-automorphic
function $G$ with $\Q$-rational Fourier coefficients, and such
that $w_N(G) = c_N/G$. To see that this gives the theorem as we
have stated it, apply the Exercise in \cite[$\S$ 5.3]{TGT}.
\\ \\
\noindent Thus, if $c_N$ is not a norm from $\Q(\sqrt{p^*})$ we
deduce that $C(N,p)$ has at least \emph{two} deficient places.  In
particular we get the following sharpening of Proposition
\ref{17}.
\begin{cor}
\label{sharp17} The deficient places of $C(17,p)$ are precisely
$\ell = p$ and $\ell = 17$.
\end{cor}
\noindent Proof: By Theorem \ref{obstruction}, $C(17,p)$ maps to
the genus zero curve $V$ with corresponding quaternion algebra
$\langle 17^3, \ p^* \rangle \cong \langle 17, \ p^* \rangle$.
Since $V(\Ql) = \emptyset$ for $\ell = p$ and $17$, \emph{a
fortiori} the same holds for $C(17,p)$. That there are no other
deficient places follows from Propositions \ref{real} and
\ref{smooth}.
\begin{thm}
\label{AtN} Suppose $N$ is prime.  Then $C(N,p)(\Q_N) = \emptyset
\iff N \equiv 1 \pmod 4$.
\end{thm}
\noindent Proof: Because $\Q(\sqrt{-2})$ and $\Q(\sqrt{-3})$ have
class number one, we may assume that $N \geq 5$. We shall apply
the work of Mazur and Rapoport on the structure of the minimal
regular model for $X_0(N)_{/\Q_N}$ \cite{MR} (especially relevant
is the picture on \cite[p. 177]{MR}). Recall that the special
fiber $M_0(N)_{/\F_p}$ of the coarse moduli space consists of two
rational curves intersecting tranversely along the supersingular
points -- with each supersingular point on
$X(1)_{/\overline{\F_N}}$ getting glued to its Galois conjugate
under the quadratic Frobenius map $\sigma: \F_{N^2} \ra \F_{N^2}$.
The Atkin-Lehner involution $w_N$ has the effect of interchanging
the two branches, and the assumption that $p^*$ is a quadratic
nonresidue modulo $N$ -- note that $(\frac{N}{p}) =
(\frac{p^*}{N})$ -- means that $C(N,p)$ is the generic fiber of an
arithmetic surface $M(N,p)_{/\Z_N}$ which is the twist of
$M_0(N)_{/\Z_N}$ by $\sigma$.  In particular, the only
$\F_N$-rational points on the special fiber of $M(N,p)$ are the
supersingular points, which are singular. \\ \indent However, as
alluded to above, $M_0(N)_{/\Z_N}$ is not necessarily a
\emph{regular} model of $X_0(N)_{/\Q_N}$.  More precisely, if $N
\equiv -1 \pmod 3$, then $j = 0$ must be blown up twice leading to
a chain of two rational curves; and similarly if $N \equiv -1
\pmod 4$, $j = 1728$ must be blown up once leading to a single
rational curve.  Now a similar procedure can be performed on
$M(N,p)_{\Z_N}$, to get a(n in fact minimal)
regular model of $C(N,p)_{/\Q_N}$,
which we shall denote by $\mathcal{C}(N,p)_{/\Z_N}$;
we need only keep track of the
effect of the twisted $\Gal(\F_{N^2}/\F_N)$ action on these
rational curves. In the $j = 0$ case, the Galois action
interchanges the two rational curves, hence leads to no new
$\F_N$-rational points. However, in the $j = 1728$ case the Galois
action the unique rational curve is evidently stabilized by the
Galois action, yielding a smooth $\F_N$-rational curve of genus
zero.  But, as above,  it is well known that every smooth genus
zero curve over a finite field is isomorphic to the projective
line, giving $N+1$ rational points on the special fiber of $\F_N$.
In summary, we have that a (minimal) regular model for $C(N,p)$ over $\Z_N$
has $\F_N$-rational points if and only if $N \equiv -1 \pmod
4$. We are done by Hensel's Lemma.
\\ \\
Remark: The proof still goes through when $N = 2$ or $N = 3$.  In
each case there is a unique supersingular point on the special
fiber.  The desingularization performed by successive blowups
of this point leads to a chain of $11$ rational curves when $N = 2$ and $5$
rational curves when $N = 3$.  Since these numbers are odd, the
Galois action stabilizes the middle element of the chain.
\\ \indent
Moreover, we have assumed $N$ to be prime only for simplicity of
exposition: for any squarefree $N$ and $\ell$ dividing $N$, the
argument gives a necessary and sufficient condition for
$C(N,p)(\Ql) = \emptyset$, namely that $(\frac{p^*}{\ell}) = -1$
and that there does not exist a supersingular point on
$X_0(N/\ell)_{/\F_{\ell}}$ whose automorphism group is divisible
by $4$.  We leave the task of converting this into an explicit
congruence condition to the interested reader.
\\ \\
\noindent Theorem \ref{AtN} gives a third proof of Proposition
\ref{17}. The relationship between Theorems \ref{obstruction} and
\ref{AtN} (whose proofs seem very different) is interesting:
neither encompasses the other, although there is a substantial
overlap: the implication $\Longleftarrow$ of Theorem \ref{AtN}
also follows from Theorem \ref{obstruction}.
\\ \\
At present I do not know whether there is a similarly simple
necessary and sufficient condition for $C(N,p)(\Q_p) = \emptyset$;
whether there is ever a deficient prime $\ell$ not dividing $Np$;
or whether $C(N,p)$ can have no deficient places but still fail to
have $\Q$-rational points.
\\ \\
Acknowledgements: I thank Dr. Nick Rogers for performing
additional calculations and Professor Shih for making me aware of
the connections with his work.  As usual, I am grateful to Noam
Elkies for rapid and helpful replies to my questions.

\end{document}